\documentclass[12pt]{article}
\usepackage[margin=2.5cm]{geometry}
\usepackage{times}
\usepackage{setspace}
\usepackage{caption}
\usepackage{lipsum} % For placeholder text

\usepackage[dvipsnames,table]{xcolor}
\definecolor{rub-blau}{RGB}{0,53,96}
\definecolor{rub-grun}{RGB}{141, 174, 16}
\definecolor{fast-rub-blau}{RGB}{160,173,189}
\usepackage{lastpage}

\usepackage{amsmath,amsthm,amssymb}
\usepackage{graphicx}
\usepackage{caption}
\usepackage{enumitem}
\usepackage[sc]{mathpazo}

\usepackage[T1]{fontenc}

\usepackage[colorlinks=true,
            linkcolor=blue,
            citecolor=blue,
            urlcolor=blue,
            pdfborder={0 0 0}]{hyperref}

\usepackage[capitalize,nameinlink,noabbrev]{cleveref}

\newtheorem{thm}{Theorem}

\newtheorem{conj}{Conjecture}

\newtheorem{remark}{Remark}

\usepackage{tikz}
\usetikzlibrary{positioning}

\usepackage{dsfont}

\usepackage{booktabs}
\usepackage{array}
\usepackage{tikz}
\usetikzlibrary{trees}

\usetikzlibrary{positioning, shapes.geometric, arrows.meta}

\makeatletter
\renewenvironment{abstract}{%
  \small
  \begin{center}%
    {\bfseries \abstractname\vspace{-.5em}\vspace{0pt}}%
  \end{center}%
  \quotation
}{%
  \endquotation
}
\makeatother

\newcommand{\firstlevel}[1]{\section{\textbf{#1}}}
\newcommand{\secondlevel}[1]{\subsection{\textbf{\textit{#1}}}}

\begin{document}

% Title
\begin{center}
    \textbf{How to Use Deep Learning to Identify Sufficient Conditions: A Case Study on Stanley's $e$-Positivity}
\end{center}

% Author information with superscript numbers
\begin{center}
    Farid Aliniaeifard\textsuperscript{1,2} (\href{mailto:farid@sdu.edu.cn}{farid@sdu.edu.cn}) \\
    Shu Xiao Li\textsuperscript{1,2} (\href{mailto:lishuxiao@sdu.edu.cn}{lishuxiao@sdu.edu.cn}) \\[0.5cm]
    
    \textsuperscript{1}\textit{Research Center for Mathematics and Interdisciplinary Sciences, Shandong University} \\
    \textsuperscript{2}\textit{Frontiers Science Center for Nonlinear Expectations, Ministry of Education, Qingdao, Shandong, 266237, P. R. China}
    \vspace{0.5cm}
\end{center}

%\usepackage[hyphens]{url}
%\usepackage{authblk}

%
%\author{Farid Aliniaeifard and Shu Xiao Li}
%\date{}
%\address{
%	Research Center for Mathematics and Interdisciplinary Sciences, Shandong University Frontiers Science Center for Nonlinear Expectations, Ministry of Education, Qingdao, Shandong, 266237, P. R. China}
%\email{farid@sdu.edu.cn}
%
%\author{Shu Xiao Li}
%
%\address{
%	Research Center for Mathematics and Interdisciplinary Sciences, Shandong University \\ Frontiers Science Center for Nonlinear Expectations, Ministry of Education, Qingdao, Shandong, 266237, P. R. China}
%\email{lishuxiao@sdu.edu.cn}
%\thanks{Both author were supported in part by the Provincial Nature Science Foundation of Shandong, Project No. ZR2024QA026 and the Fundamental Research Funds for the Central Universities.}
%%
%\author{Farid Aliniaeifard\thanks{Both authors were supported in part by the Provincial Nature Science Foundation of Shandong, Project No.
%ZR2024QA026 and the Fundamental Research Funds for the Central Universities. } \and Shu Xiao Li}

%\address{\addressmark{1} Research Center for Mathematics and Interdisciplinary Sciences, Shandong University  \\ \addressmark{2}Frontier Science Center for Nonlinear Expectations, Ministry of Education, Qingdao,  Shandong, 266237, P. R. China}
%

%\maketitle

\begin{abstract} In a study, published in \emph{Nature},  researchers from DeepMind and mathematicians demonstrated a general framework using machine learning to make conjectures in pure mathematics. Here, we build upon this framework to develop a method for identifying sufficient conditions that imply a given mathematical statement. As a demonstration, we apply this process to Stanley's problem of $e$-positivity of graphs--one of the problems that has been at the center of algebraic combinatorics for the past three decades.  Guided by AI, we rediscover that one sufficient condition for a graph to be $e$-positive is that it is co-triangle-free. Based on Saliency Map analysis, we suggest that the classification of $e$-positive graphs is more related to continuous graph invariants rather than the discrete ones, which we support it with three conjectures. 
Furthermore, we show that the claw-free and claw-contractible-free graphs with $10$ and $11$ vertices are $e$-positive, resolving a conjecture by Dahlberg, Foley, and van Willigenburg. 
\end{abstract} 

% Add keywords after the abstract

\firstlevel{Introduction} 
How do we make conjectures in combinatorics, or mainly in discrete mathematics? Specifically, how do we find conditions $P$ that imply a given statement $Q$? How do we know if a chosen condition is optimal, and if a weaker one exists that still implies $Q$, and even might suggest a proof strategy?  Before the age of AI, the answers lay in a combination of human intuition and exhaustive computations. Now, machine learning algorithms, particularly those having attribution techniques,  provide a new framework for generating and proving conjectures. 

Mathematicians like Marc Lackenby, Geordie Williamson, and DeepMind scientists like the Nobel Laureate Demis Hassabis proposed a general framework %(see Figure \ref{fig:Nature}) 
for guiding mathematicians' intuition to obtain impactful mathematical results \cite{Nature21}. In their work, the framework helped reveal a new connection between the algebraic and geometric structures of knots and provided a refinement of the combinatorial invariance conjecture for symmetric groups. See more breakthroughs that demonstrated AI’s capability to assist mathematicians in \cite{alfarano2024global, coolsaet2023house, douglas2022numerical, georgiev2025mathematical, he2025murmurations,   wang2023asymptotic, wang2025discovery, wu2016counting}. 

In many problems in mathematics, including positivity problems, we often know the target condition $Q$ and seek sufficient conditions $P$ that imply $Q$. An example is identifying which graphs are $e$-positivity \cite{stanley1995symmetric}\footnote{We give enough details about Stanley's chromatic symmetric functions in this footnote. Let $G$ be a graph and $\mathbb{Q}[[x_1,x_2,\dots ]]$ be the ring of power series in the countable commuting variables $\{x_1,x_2, \dots \}$. The \emph{chromatic symmetric function} of $G$ is $X_G=\sum_{\kappa} \prod_{v\in V(G)} x_{\kappa(v)}.$ This is a symmetric function. One of the well-known basis of the space of symmetric functions is the \emph{elementary}  basis $\{e_\lambda: \lambda~ \text{is a partition}\}$, where for $\lambda=(\lambda_1,\dots, \lambda_k)$, $$e_\lambda=\left(\sum_{i_1<\dots<i_{\lambda_1} }x_{i_1}\dots x_{i_{\lambda_1}}\right) \cdots \left(\sum_{i_{n-\lambda_k+1}<\dots<i_{n} }x_{i_{n-\lambda_k+1}} \dots x_{i_{n}}\right).$$ A graph is said to be \emph{$e$-positive} if its chromatic symmetric function can be written as a nonnegative linear combination of the elements of the elementary basis.}. A complete classification of $e$-positive graphs is a very difficult problem and still open. Even for the family of incomparability graphs of $(3+1)$-free posets, conjectured by Stanley and Stembridge in 1993 \cite{stanley1993immanants}, after three decades in 2024, Hikita \cite{hikita2024proof} showed the conjecture is true.

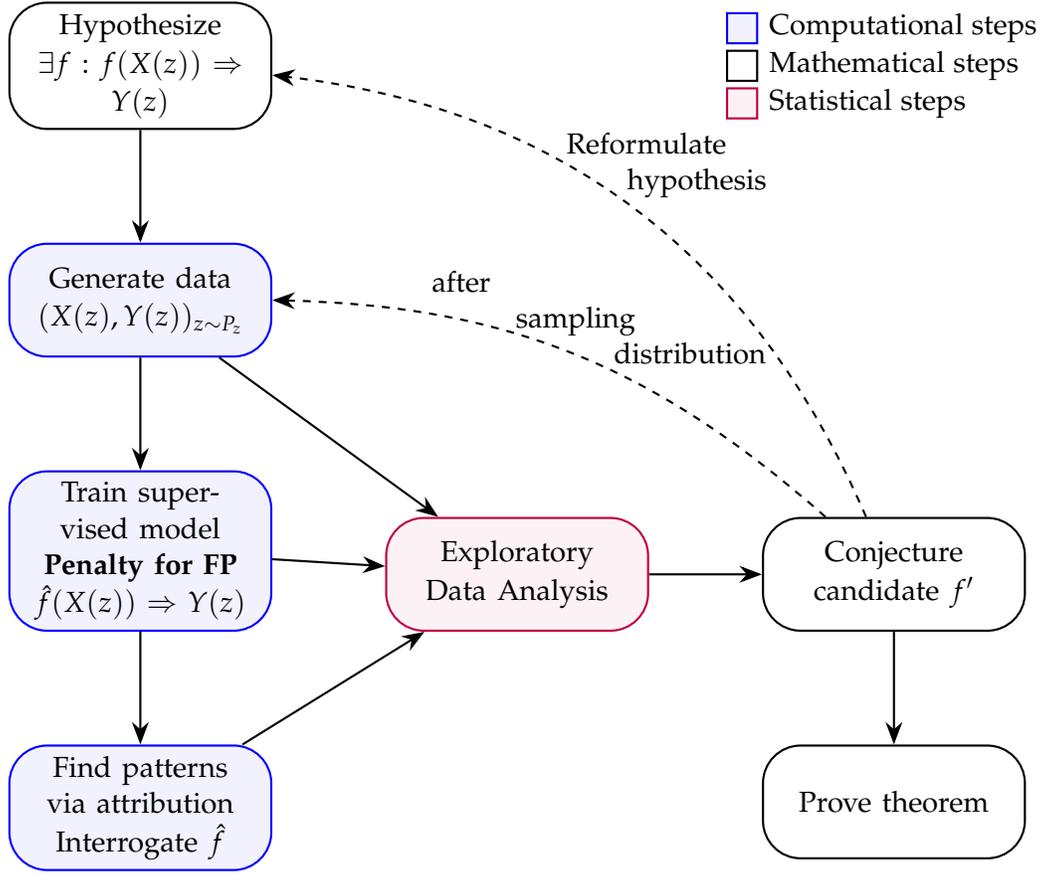
\begin{figure}[t]
    \centering
\begin{tikzpicture}[
    node distance=1.5cm and 1.5cm,
    stepnode/.style={
        rectangle,
        rounded corners=5mm,
        minimum width=3cm,
        minimum height=1.5cm,
        text width=3.2cm,
        align=center,
        draw=blue,
        thick,
        fill=blue!5
    },
     statnode/.style={
        rectangle,
        rounded corners=5mm,
        minimum width=3cm,
        minimum height=1.5cm,
        text width=3.2cm,
        align=center,
        draw=purple,
        thick,
        fill=purple!5
    },
    mathnode/.style={
        rectangle,
        rounded corners=5mm,
        minimum width=3cm,
        minimum height=1.5cm,
        text width=3.2cm,
        align=center,
        draw=black,
        thick,
        fill=white
    },
    dashedarrow/.style={
        -{Stealth[length=3mm]},
        dashed,
        thick
    },
    solidarrow/.style={
        -{Stealth[length=3mm]},
        thick
    }
]

% Nodes
\node[mathnode] (A) {\small Hypothesize\\ $\exists f: f(X(z))\Rightarrow Y(z)$};
\node[stepnode, below=of A] (B) {\small Generate data $(X(z),Y(z))_{z\sim P_z}$};
\node[stepnode, below=of B] (C) {\small Train supervised model \\ {\bf Penalty for FP} \\ $\hat{f}(X(z)) \Rightarrow Y(z)$};
\node[stepnode, below=of C] (D) {\small Find patterns via attribution\\Interrogate $\hat{f}$};
\node[statnode, above right=of D] (E) {\small Exploratory Data Analysis};
\node[mathnode, right=of E] (F) {\small Conjecture candidate $f'$};
\node[mathnode, below=of F] (G) {\small Prove theorem};

% Main flow arrows
\draw[solidarrow] (A) -- (B);
\draw[solidarrow] (B) -- (C);
\draw[solidarrow] (C) -- (D);
\draw[solidarrow] (D) -- (E);
\draw[solidarrow] (E) -- (F);
\draw[solidarrow] (F) -- (G);

% Additional arrows
\draw[solidarrow] (B) -- (E);
\draw[solidarrow] (C) -- (E);
%\draw[solidarrow] (D) -- (E);

% Dashed arrows with labels
\draw[dashedarrow] (F) to[bend right=20] node[midway,right,align=left,font=\small] 
   { \hspace{-2cm}after \\ \hspace{-2cm} ~~~~~~~~~sampling\\~ \hspace{0.2cm}distribution~~ \\ ~\\ ~} (B);
\draw[dashedarrow] (F) to[bend right=30] node[midway,above,align=center,font=\small] 
    {~~~Reformulate\\ ~~~~~~~~~~~~~~~hypothesis} (A);

% Legend
\node[rectangle, draw=black, thick, fill=white, minimum size=0.4cm, label=right:{\small Mathematical steps}] 
    (mathlegend) at (8,0) {};
\node[rectangle, draw=blue, thick, fill=blue!5, minimum size=0.4cm, label=right:{\small Computational steps}] 
    (complegend) at (8,0.5) {};
\node[rectangle, draw=purple, thick, fill=purple!5, minimum size=0.4cm, label=right:{\small Statistical steps}] 
    (complegend) at (8,-0.5) {};
\end{tikzpicture}
    \caption{The framework for finding sufficient conditions}
    \label{fig: sc}
\end{figure}

In this study, we present a general framework (see Figure \ref{fig: sc}) that uses machine learning and attribution techniques to guide mathematicians' intuition in finding sufficient conditions $P$ that imply a statement $Q$. This is achieved by modifying the classical loss functions. The core idea is that, given $Q$, we create a dataset of many possible conditions that might influence it. We then design deep learning models that impose a severe penalty on false positive predictions, with the goal of achieving a model with $100\%$ precision. Then we use Saliency Map analysis to identify conditions $Q$. Using this framework in the study of $e$-positive graphs, we observe that the $e$-positivity of graphs is more closely related to continuous graph invariants. Based on this, in Section \ref{sec:conjectures}, we present several conjectures that connect $e$-positivity to continuous invariants.

In Section \ref{sec: AI}, we present our general framework for finding sufficient conditions that imply $Q$. In Section \ref{sec: e-pos}, we use our framework and rediscover that the co-triangle-free graphs are $e$-positive, and any graph with $n$ vertices and independence number $\lceil \frac{n}{2}\rceil +1$ is not  $e$-positive. In addition to that, using Saliency Map analysis of neural networks, we show that a complete classification of $e$-positive graphs depends more on continuous graph invariants, supported by several conjectures in Section \ref{sec:conjectures}.  We also provide data confirming that the claw-free and claw-contractible-free graphs with $10$ and $11$ vertices are $e$-positive, resolving the conjecture in \cite[page 22]{dahlberg2020resolving}. In Section \ref{method}, we detail the methodology used to obtain these results. In the last section, we discuss the advantages of our framework, the choice of the $e$-positivity problem, and, furthermore, what Saliency Map analysis reveals about the relationship between continuous variables and $e$-positivity.

\firstlevel{Using AI to Identify Sufficient Conditions}\label{sec: AI}
We first establish our terminology and notation. 
A \emph{property} of a category $\mathcal{C}$ is a categorical invariant, that is, an assignment from $\mathcal{C}$ to a set $I$, $$X: \mathcal{C}\rightarrow I,$$ such that for objects $z$ and $z'$ in $\mathcal{C}$ with $z\cong z'$, then $X(z)=X(z')$. 
For example, consider the category of graphs $\mathcal{G}$. The assignment
$${\rm inc}_{(3+1)}(G): \mathcal{G} \rightarrow \{0,1\},$$ defined by  ${\rm inc}_{(3+1)}(G)=1$ if $G$ is the incomparability graph of a $(3+1)$-free poset and $0$ otherwise, is a property of graphs.   Mathematicians often make conjectures that predict some relation between properties $X(z)$ and $Y(z)$ of the objects $z$, typically of the form: 
$$\forall z \in \mathcal{C} \text{~s.t~} X(z) \in A \Rightarrow Y(z)\in B.$$ 
For example, Stanley and Stembridge  \cite{stanley1993immanants} made the following conjecture in 1993: 
 $$\forall G\in \mathcal{G} \text{~s.t~} {\rm inc}_{(3+1)}(G)\in \{1\} \Rightarrow  e\_{\rm pos}(G) \in  \{1\},$$
 where  $e\_{\rm pos}$ is an assignment whose value is $1$ for $e$-positive graphs and $0$ otherwise. 
 
 Sometimes, we have a question in mathematics that states for which objects $z$, a property $Y(z)\in B$ holds. That is, can we find sufficient conditions of the form $X(z)\in A$ such that
 \begin{equation}\label{conj}
 \forall z\in\mathcal{C} \text{~ s.t~}X(z)\in A \Rightarrow Y(z)\in B.
 \end{equation}
Examples of such problems are positivity problems, which ask, given an element $v$ in a vector space with a specific basis, whether the coefficients in the expansion of $v$ in that basis are non-negative.  More specifically, for chromatic symmetric functions, the question is to find a graph property $X$ and  a set $A$ such that 
$$\forall G\in\mathcal{G} \text{~s.t}~ X(G)\in A \Rightarrow e\_{\rm pos}\in \{1\}.$$
 
To present the general framework for finding sufficient conditions, we must first understand the metrics used for evaluations of predictions, as any conjecture is a form of prediction. If we want to make a conjecture of the form \eqref{conj}, which metric should we prioritize?  

\secondlevel{Metrics and conjectures}  Any conjecture is a prediction, and specific metrics provide a comprehensive understanding of its quality. For example, consider a test that predicts whether a person has cancer. If the test predicts cancer and the person indeed has cancer, then the prediction is called a \emph{true positive (TP)}. Let $N_{\text{TP}}$ be the number of true positives. If the test predicts cancer but the person does not have it, this is called \emph{false positive (FP)}. Let $N_{\text{FP}}$ be the number of false positives. Similarly, we can define \emph{true negative (TN)}, $N_{\text{TN}}$ and \emph{false negative (FN)}, $N_{\text{FN}}$. 
\begin{table}[ht]
\centering
\begin{tabular}{|c|l|}
\hline 
${\text{precision}} = {N_{\text{TP}}}/({N_{\text{TP}}+ N_{FP}})$ & ${\text{recall}} = {N_{\text{TP}}}/({N_{\text{TP}} + N_{\text{FN}}})$\\
\hline 
${\text{accuracy}} = (N_{\text{TP}}+N_{\text{TN}}) / Total $& ${\text F1\text{-Score}} = 2  (\text{percision} \times \text{recall}) / (\text{percision} + \text{recall})$\\
\hline 
\end{tabular} 
\captionof{table}{Metrics for a comprehensive understanding of a binary prediction}
\label{tab:metrics}
\end{table}

Our goal in the next section is to set up a deep learning algorithm that, when we know $Y(z)\in B$, suggests a conjecture of the form \eqref{conj}.
Let's check the metrics in Table \ref{tab:metrics} to see which one guarantees that the conjecture holds.  In this context, we can only derive the following information from a conjecture of the form \eqref{conj}, 
$$
\begin{tabular}{|c|l|}
\hline 
$N_{TP}$ & The number of objects $z$, that $X(z) \in A$ and $Y(z)\in B$. \\
\hline
$N_{FP}$ & The number of objects $z$, that $X(z) \in A$ and $Y(z)\not\in B$. \\
\hline
%$N_{TN}$ & The number of objects $z$, that $X(z) \not\in A$ and $Y(z)\not\in B$. \\
%\hline
%$N_{FN}$ & The number of objects $z$, that $X(z) \not\in A$ and $Y(z)\in B$. \\
%\hline
\end{tabular}
$$
Note that such a conjecture, as a prediction, gives partial information about $N_{TN}$ and $N_{FN}$. Therefore, we can say that a conjecture of the form \eqref{conj} is a prediction that aims for $100\%$ percision when checking if $Y(z)\in B$, i.e., $N_{FP}=0$.  

\begin{remark} 
Consider that, additionally, a conjecture of the form 
$$X(z)\in A \Leftrightarrow y(z)\in B,$$ predicts that $N_{FP}=0$ and $N_{FN}=0$, meaning all its metrics in Table \ref{tab:metrics} are $100\%$. 
\end{remark} 

 Thus, when Stanley and Stembridge made the following conjecture: 
 $$\forall G\in\mathcal{G} \text{~s.t~} {\rm inc}_{(3+1)}(G)\in \{1\} \Rightarrow  e\_{\rm pos}(G) \in  \{1\},$$
 they were making a prediction with $100\%$ precision.

From our discussion in this section, if we want to set up a deep learning algorithm to find sufficient conditions for a statement, we need a model that makes predictions with $100\%$ precision. 
For any deep learning model, one can design a loss function that prioritizes performance on a specific metric.  

\secondlevel{The general framework}
In this section, we outline a general framework for identifying sufficient conditions for $Y(z)\in B$ using a set of properties. 
The main differences between our framework and the general framework in \cite{Nature21} are the inclusion of a severe penalty for FP predictions and also an emphasis on exploratory data analysis. 

The general framework proceeds as follows:
\begin{enumerate} 
 \item Propose the hypothesis that the properties $$X(z) = X_1(z), X_2(z),\cdots, X_n(z)$$ can help with identifying the target property $Y(z)\in B$.  
 
\item Generated a dataset $X_1(z),X_2(z),\cdots, X_n(z), Y(z)$ for a large, representative family of the objects $z$.

\item Train a neural network (or another suitable supervised machine learning model) to predict $Y(z)$ using the features $X_1(z), X_2(z),\dots, X_n(z)$ as input. The model is explicitly biased towards high precision by employing a loss function that assigns a severe penalty to FP predictions. 

The model outputs a prediction $\hat{f}(X_1,X_2,\cdots, X_n)$.

\item If the model achieves a prediction greater than a random baseline, even if marginally, while maintaining precision close to $1$, this provides evidence that the features $X(z)$ could help with identifying $Y(z)\in B$. 

\item Perform Saliency Map analysis to rank the features $X_i$ based on their contribution to the model's prediction, identifying a subset $X_{i_1}, X_{i_2}, \dots, X_{i_q}$ with the highest importance. The Saliency Map analysis \cite{simonyan2013deep} for an input feature $X_i$ and dataset $X$ computes $$s_i = \frac{1}{|\mathcal{X}|}\sum_{x\in \mathcal{X}}\left|\frac{\partial\hat{ f}}{\partial X_i} (x)\right|$$ and then rank features using $s_i$. 

\item Conduct exploratory data analysis on the high-importance features $X_{i_1}, X_{i_2}, \dots, X_{i_q}$ to identify specific sets $S_1,S_2,\dots, S_q$ such that the following conjecture $f'$ can be made: 
 $$\forall z \in \mathcal{C} \text{~s.t~}X_{i_1}(z)\in S_1, X_{i_2}(z)\in S_2, \dots, X_{i_q}(z)\in S_q \Rightarrow Y(z)\in B.$$
 
 \item The process may be repeated several times before settling on a viable conjecture. The selection of conjectures should be both engaging and reasonably believable; here, machines are less capable, and mathematicians play an important role.

\end{enumerate} 

\firstlevel{Stanley's $e$-positivity}\label{sec: e-pos}
Stanley and Stembridge in \cite{stanley1993immanants} conjectured that the incomparability graphs of $(3+1)$-free posets are $e$-positive. Guay-Paquet \cite{guaypaquet2013modular} later showed the problem could be reduced to unit interval graphs, and Hikita \cite{hikita2024proof} recently presented a proof of this conjecture. This problem has been one of the central problems in algebraic combinatorics for three decades. Using AI, we identified a potentially larger family of $e$-positive graphs.

\begin{figure}[h]
    \centering
    \includegraphics[width=0.7\textwidth]{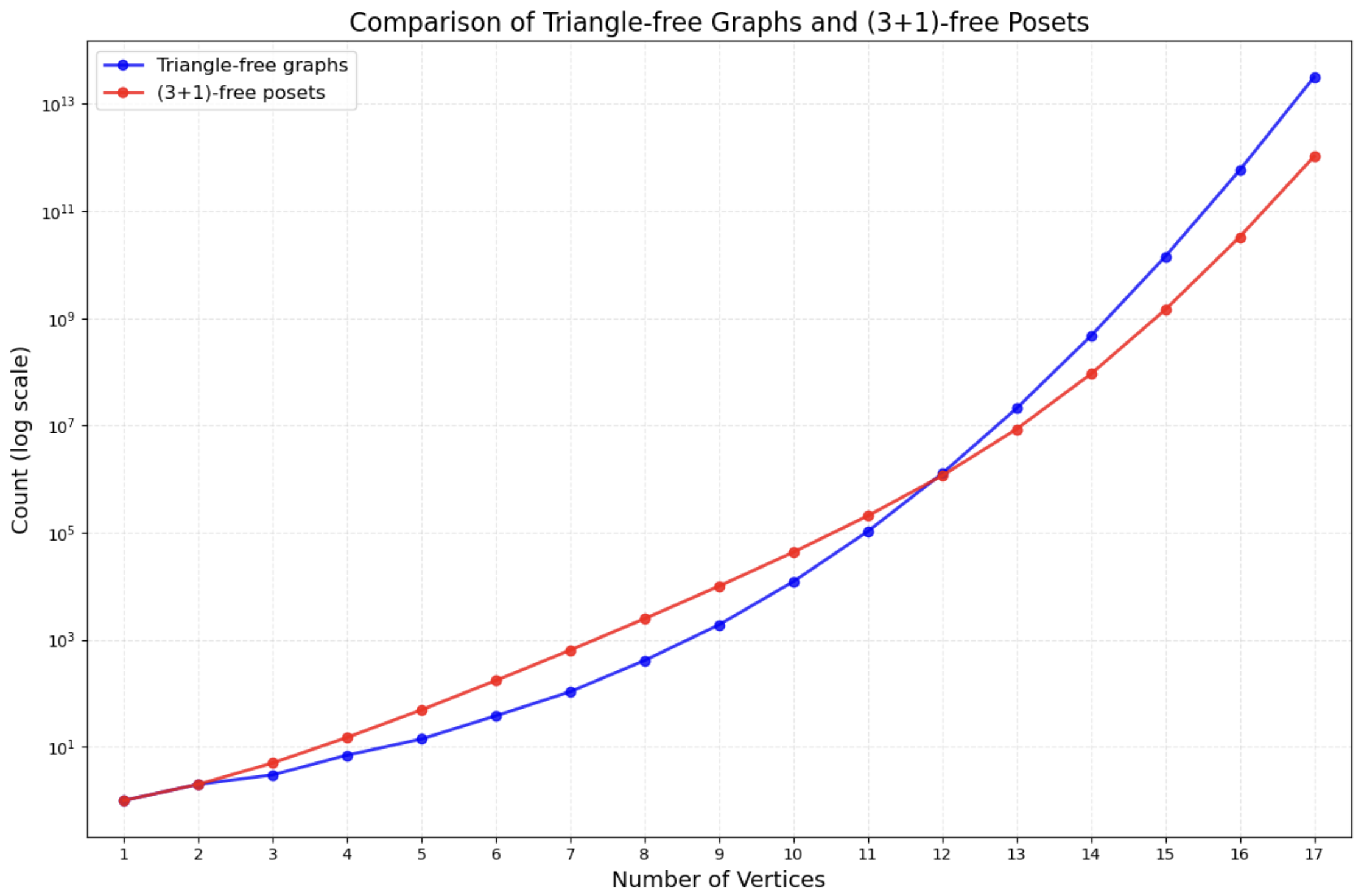}
     \caption{Comparing the number of triangle-free graphs and $(3+1)$-free posets}
     \label{fig:ss5}
\end{figure}

Our hypothesis was that there are more connections between graph invariants and $ e$-positivity. A supervised machine learning model detected graphs whose complements are triangle-free are $e$-positive. This is significant, as shown in Figure \ref{fig:ss5}, the number of triangle-free graphs (OIES A006785) far exceeds the number of graphs considered in the $(3+1)$-free conjecture (see \cite{guaypaquet2013structure}).

\begin{figure}[h]
    \centering
    \begin{minipage}{0.48\textwidth}
        \centering
        \includegraphics[width=\textwidth]{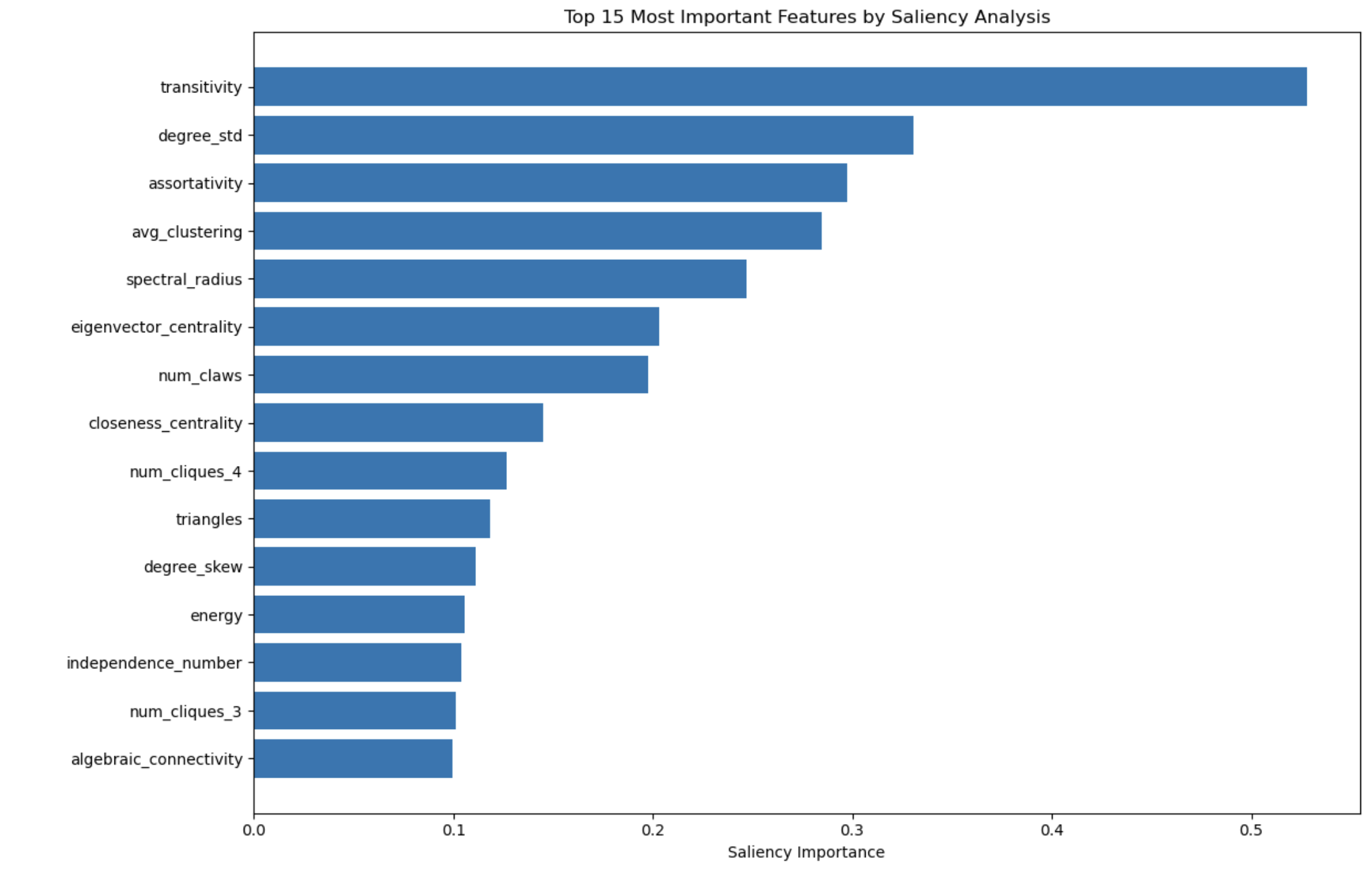}
        %\caption{Saliency map visualization}
        %\label{fig:ss_a}
    \end{minipage}
    \hfill
    \begin{minipage}{0.48\textwidth}
        \centering
        \includegraphics[width=\textwidth]{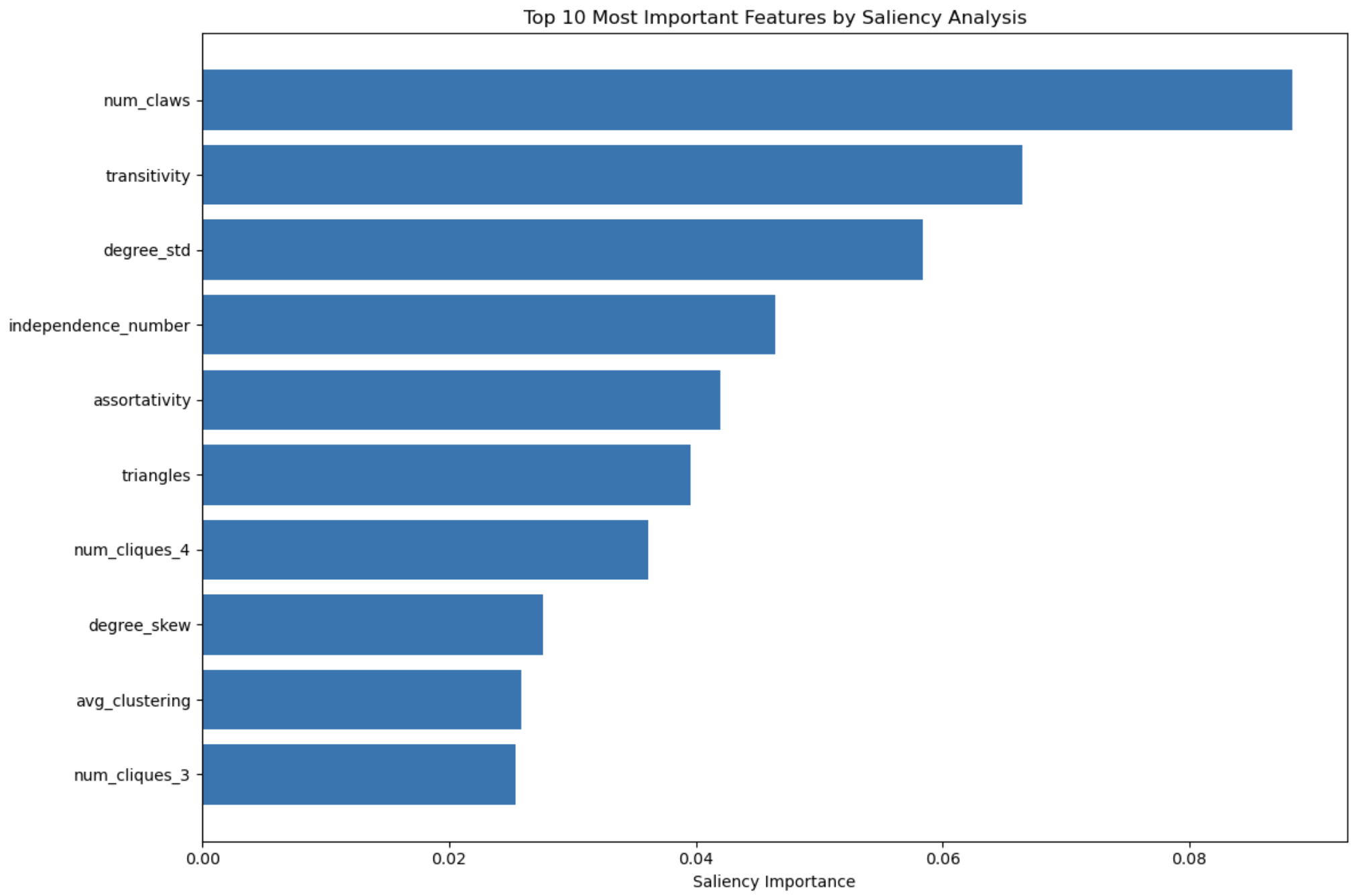}
        %\caption{Saliency map visualization}
        %\label{fig:ss_b}
    \end{minipage}
    \caption{Saliency Map visualizations}
    \label{fig:ss}
\end{figure}

We first trained a model on numerous graph invariants, which achieved 94.3\% accuracy in predicting whether a graph is $e$-positive.
Using attribution techniques, Saliency Map analysis, we identified the fifteen most relevant features for this binary classification; these are shown in Figure \ref{fig:ss}.
%\begin{figure}[t]
%    \centering
%    \includegraphics[width=0.7\textwidth]{Saliency1}
%    \caption{Saliency map visualization}
%    \label{fig:ss}
%\end{figure}
One interesting observation here is that the six most important graph invariants for the model to determine if a graph is $e$-positive are continuous. This is evidence that $e$-positivity of graphs is more related to continuous graph invariants. Look at Section \ref{sec:conjectures} for supporting conjectures. 

\newpage 

We then trained a model on these fifteen graph invariants, this time focusing on precision, and setting the model to severely penalize FP. This approach increases the likelihood of identifying a conjecture of the form \eqref{conj}, that is, if certain graph invariants fall within specific ranges, the graph is $e$-positive. 

One of the interesting observations is that the model identifies the number of claws as the most important graph invariant in determining $e$-positivity. This is fascinating because the Stanley-Stembridge conjecture concerned certain claw-free graphs, and also Stanley shows that the co-triangle-free graphs are $e$-positive; here, the intuition of mathematicians matches the computations of deep learning models three decades later. \\

\begin{figure}[h]
    \centering
    \includegraphics[width=0.8\textwidth]{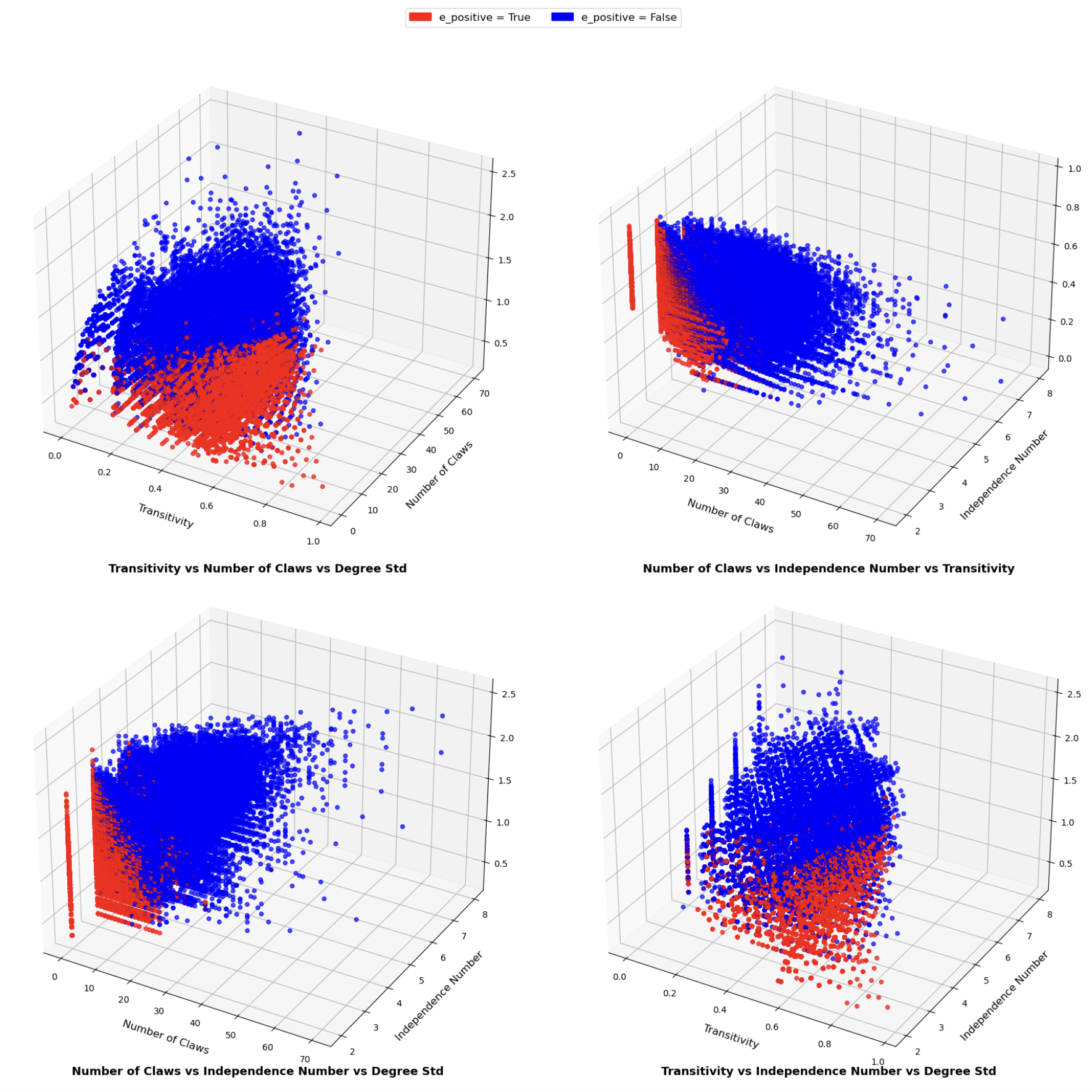}
    \caption{$e$-positivity and the four most important graph invariants}
    \label{fig:Final}
\end{figure}

The formal philosophy of Exploratory Data Analysis (EDA) was pioneered by John Tukey in the 1970s \cite{tukey1977exploratory}. EDA involves the visualization of the data to understand its features and discover relations between them. After statistically looking at data and visualizing the four most important features, number of claws, transitivity, degree standard deviation, and independence number, as you see in Figure \ref{fig:Final}, we could make two conjectures.
The first conjecture was that if a graph has independence number $2$--or equivalently, if its complement is triangle-free--then the graph is $e$-positive. The second conjecture was that if a graph with $n$ vertices has an independence number of $\lceil \frac{n}{2}\rceil+1$, then it is not $e$-positive. We then searched the literature and found the following results. For a reference for Theorem \ref{thm: tirangle-free}, see \cite[Table 1]{banaian2024epositivity}, which presents a list of $e$-positive graphs, and Theorem \ref{thm: e-neg} can be proven by showing that such a graph has no connected partition of the shape $\lambda$ where all parts of $\lambda$ are $2$ or except the last part all parts are $2$, and so by  \cite[Proposition 1.3.3]{Wolfgang1997} is not $e$-positive. 
 
\begin{thm}\label{thm: tirangle-free}
If the complement of a graph is triangle-free, then the graph is $e$-positive. 
\end{thm}

\begin{thm}\label{thm: e-neg}
If the independence number of a graph with $n$ vertices is $\lceil \frac{n}{2}\rceil+1$, then it is not $e$-positive. 
\end{thm}

When we figured out that the complements of triangle-free graphs are $e$-positive, we found that these graphs are claw-free and claw-contractible-free. Furthermore, EDA of claw-free and claw-contractible-free graphs revealed that most of them have small independence numbers. We therefore developed a program optimized for such graphs, computed the remaining cases, and showed that they are $e$-positive. The chromatic symmetric functions of all graphs with $10$ and $11$ vertices that are claw-free and claw-contractible free are available at \href{https://github.com/AIMath-Lab/SufficientConditions}{https://github.com/AIMath-Lab/SufficientConditions}.

\firstlevel{Methodology}\label{method}

This study employed a two-stage machine learning methodology to identify the graph invariants determining $e$-positivity. The overall process consisted of an initial feature selection phase using a high-accuracy model, followed by a model focused on maximizing precision.

\secondlevel{Data generation and initial feature set}

Using SageMath and the NetworkX Python package, a dataset of all connected, non-isomorphic graphs with $9$ vertices was generated. A binary label for $e$-positivity was determined, and 44 graph invariants were computed. These graph invariants contain the topological, spectral, and combinatorial properties.

\secondlevel{Stage 1: Comprehensive model for feature importance analysis}
The first stage aimed to identify the most important graph invariants from the initial set of $44$. The model architecture used for the experiment was a fully connected, feed-forward neural network with hidden sizes [256, 128, 64, 32] and ReLU activations, with Batch Normalization regularization, compiled with the Adam optimizer.
 The task was framed as a binary classification problem, cross-entropy loss as an optimizable loss function, and test classification accuracy as a metric of performance. 

This model achieved a high predictive accuracy of \textbf{94.3\%}, which indicates the existence of underlying patterns associated with $e$-positivity. To find the most important graph-invariants impacting $e$-positivity, we performed a Saliency Map analysis, which computes the average of the gradient of the model's output $\hat{f}$ with respect to its input features, that is for an input feature $X_i$ and dataset $\mathcal{X}$, it computes $$s_i = \frac{1}{|\mathcal{X}|}\sum_{x\in \mathcal{X}}\left|\frac{\partial\hat{ f}}{\partial X_i}(x) \right|.$$ We then rank features using $s_i$. From this analysis, the \textbf{top 15 most important features}, shown in Figure \ref{fig:ss} were identified for the next stage of modeling.

\secondlevel{Stage 2: Precision-optimized model on selected features}
The objective of the second stage was to build a model that guarantees high reliability for its positive predictions, as measured by precision. To this end, we trained an identical neural network architecture, except adding dropout normalizers, but using only the top 15 features identified in Stage 1.

A critical modification was made to the loss function.  A \textbf{custom weighted binary cross-entropy loss function} was designed to impose a severe penalty on FP predictions, effectively instructing the model to be highly conservative when classifying a graph as $e$-positive. The model successfully achieved the primary goal of this stage, reaching \textbf{100\% precision} on the test set. This indicates that whenever the model predicts a graph to be $e$-positive, it can be trusted with high confidence. The overall accuracy for this precision-optimized model was \textbf{66.4\%}.
Then we perform another Saliency Map analysis to identify the top four features, which are the number of claws, transitivity, degree standard deviation, and independence number. 
\cite{Nature21}

\firstlevel{Conjectures}\label{sec:conjectures}
In this section, we use several of the most important invariants derived from Saliency Map analysis to formulate conjectures. In these conjectures, we suggest how continuous invariants such as transitivity (three times the number of tringles in the graph divided by the number of connected triples), assortativity (see \cite[Equation (21)]{Newman}), and degree distribution are related to the $e$-positivity of graphs.To the best of our knowledge, this is the first time in the literature that continuous graph invariants have been related to the $e$-positivity of graphs.

Let $\mathcal{G}_n$ denote the set of all connected graphs on $n$ vertices. 

\begin{conj}[Stochastic Ordering of Assortativity by $e\_{\rm pos}$]
For graphs in $\mathcal{G}_n$, the conditional distribution of assortativity given $e\_{\rm pos} = 1$ is stochastically greater 
than the conditional distribution given $e\_{\rm pos}= 0$. 

Formally, for $G\in \mathcal{G}_n$, let $X(G)\in [-1,1]$ be the assortativity of $G$. Then
 for all $t \in [-1,1]$,
\[
\mathbb{P}\bigl( X(G) > t \mid e\_{\rm pos}(G) = 1 \bigr) \;\ge\; \mathbb{P}\bigl( X(G) > t \mid e\_{\rm pos}(G) = 0 \bigr),
\]
and there exists at least one threshold $t_{0} \in [0,1]$ for which the inequality is strict. 
Consequently, the mean assortativity satisfies
\[
\mathbb{E}\bigl[ X(G) \mid e\_{\rm pos}(G) = 1 \bigr] \;>\; \mathbb{E}\bigl[X(G) \mid e\_{\rm pos}(G) = 0 \bigr].
\]
\end{conj}

\begin{figure}[h]
    \centering
     \includegraphics[width=0.48\textwidth]{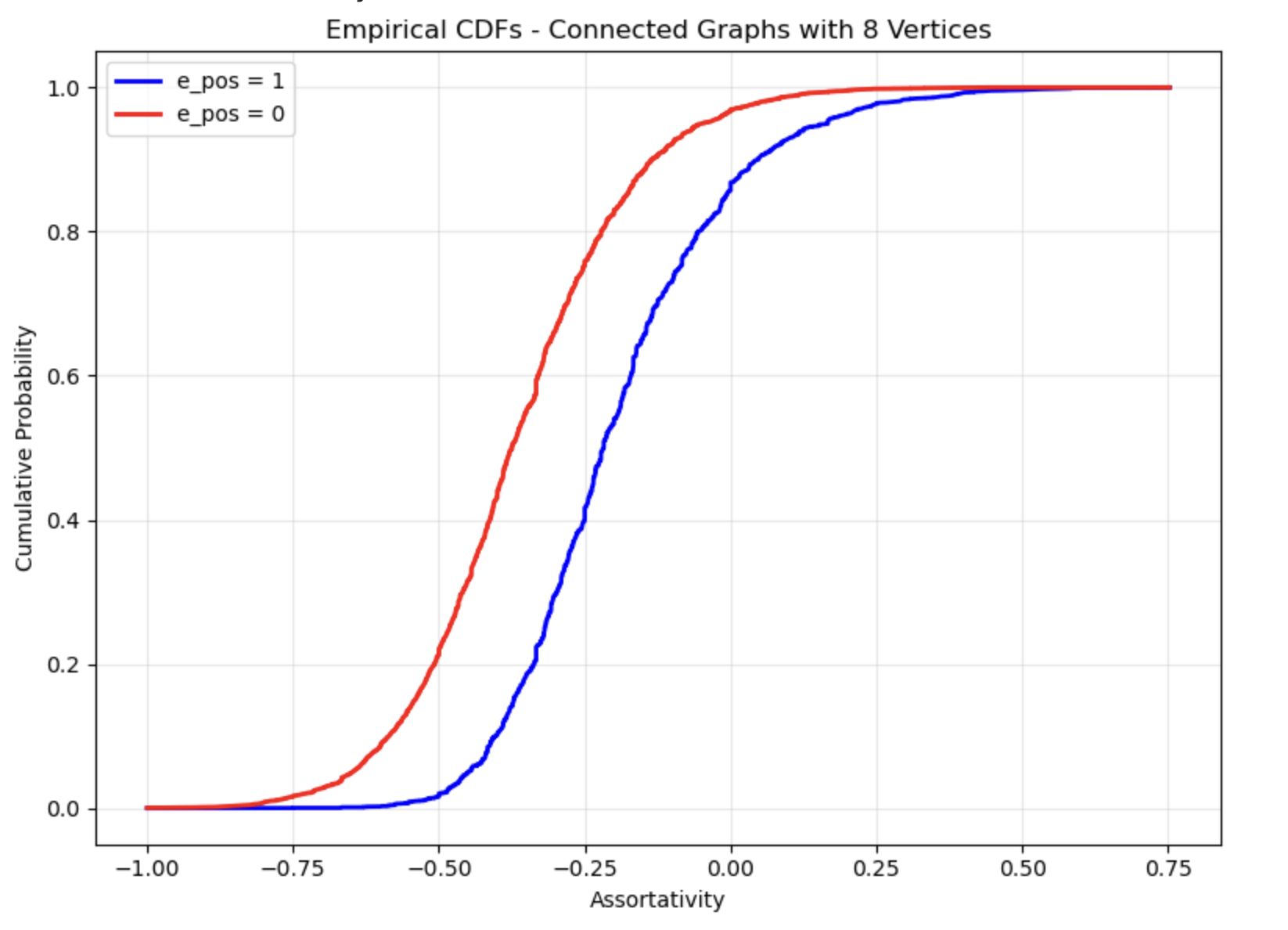} ~~ \includegraphics[width=0.48\textwidth]{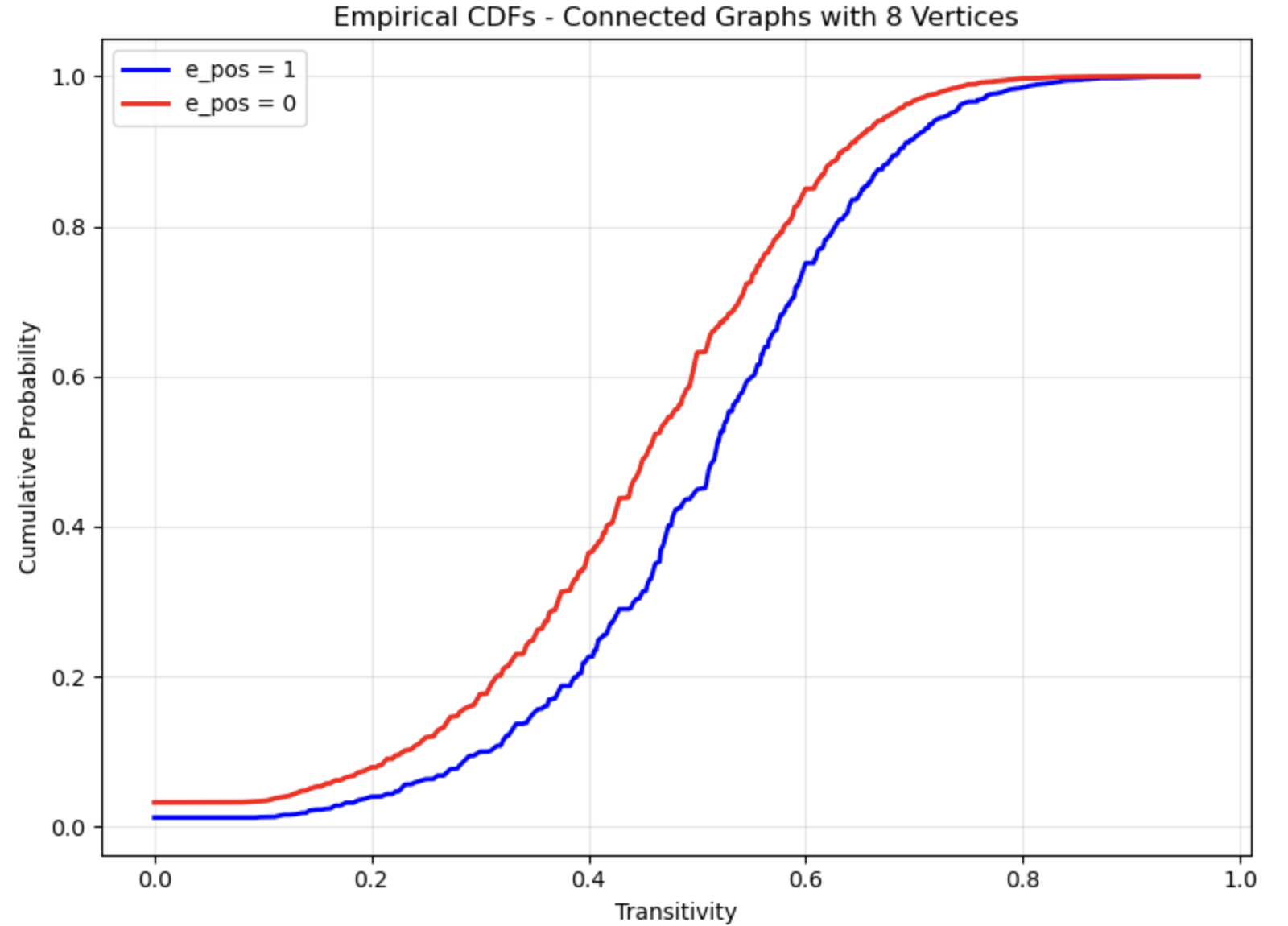} \\
        \includegraphics[width=0.48\textwidth]{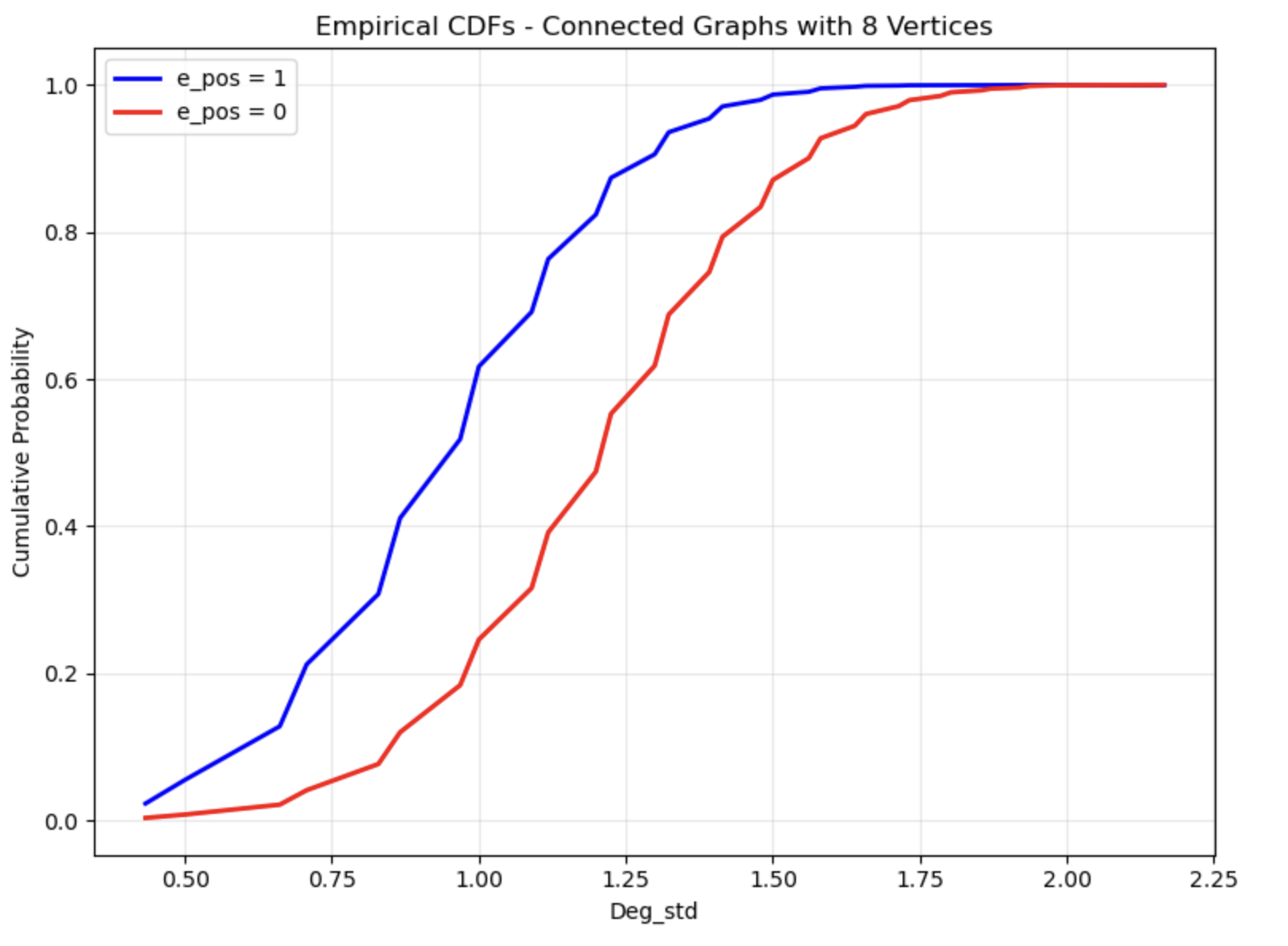}  ~~ \includegraphics[width=0.48\textwidth]{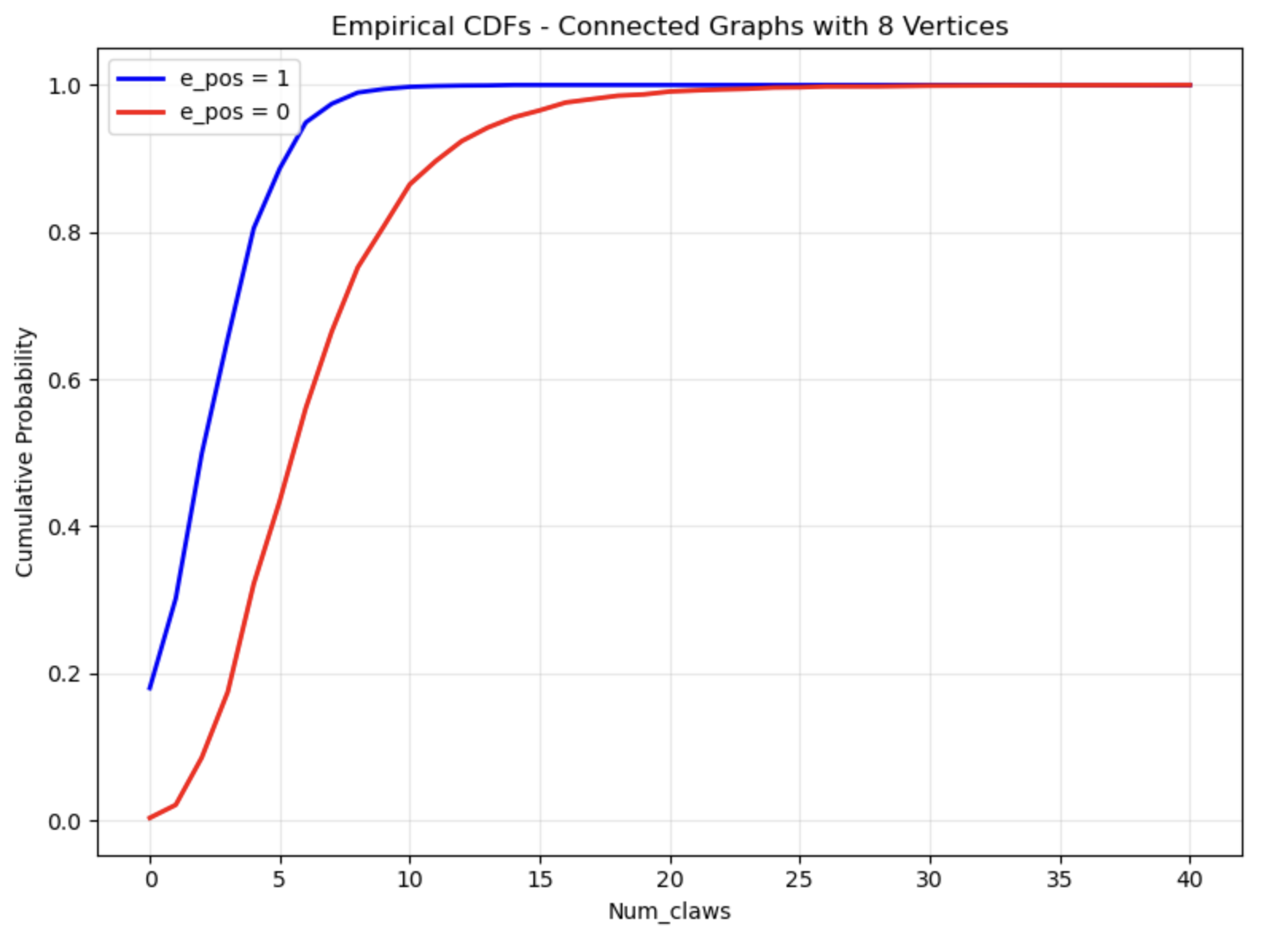}
    \caption{Empirical CDFs showing the distribution of assortativity, transitivity, degree distribution, and number of claws for $e\_{\rm pos}=1$ (blue) and $e\_{\rm pos}=0$ (red)}
    \label{fig:stochasticDominance}
\end{figure}

The above conjecture, together with Bayes’ theorem,
$$
\mathbb{P}(e_{\rm pos}=1 \,|\, X(G)>t)
= \frac{\mathbb{P}(X(G)>t \,|\, e\_{\rm pos}=1)\, \mathbb{P}(e\_{\rm pos}=1)}{\mathbb{P}(X(G)>t)},
$$
and the fact that
$$
1 = \mathbb{P}(e\_{\rm pos}(G)=1) + \mathbb{P}(e\_{\rm pos}(G)=0),
$$
implies that
$$
\mathbb{P}(e\_{\rm pos}=1 \,|\, X(G)>t) \geq \mathbb{P}(e\_{\rm pos}=1).
$$

Moreover, for some value of $t$ we have $\mathbb{P}(e\_{\rm pos}=1 \,|\, X(G)>t) = 1$. It is natural to ask what this threshold is.
Note that if the assortativity of a graph is bigger than the threshold, then the graph is $e$-positive.

%\begin{figure}[h]
%    \centering
%     \includegraphics[width=0.6\textwidth]{ass}
%    \caption{ }
%    \label{fig:}
%\end{figure}

\begin{conj}[Stochastic Ordering of Transitivity by $e\_{\rm pos}$]
For graphs in $\mathcal{G}_n, n\geq 6$, the conditional distribution of transitivity given $e\_{\rm pos}(G) = 1$ is stochastically greater 
than the conditional distribution given $e\_{\rm pos}(G)= 0$.
%\begin{figure}[h]
%    \centering
%    \includegraphics[width=0.6\textwidth]{trans}
%    \caption{ }
%    \label{fig:}
%\end{figure}
\end{conj}

\begin{conj}[Stochastic Ordering of Degree Distribution by $e\_{\rm pos}$]
For graphs in $\mathcal{G}_n$, the conditional distribution of degree distribution given $e\_{\rm pos}(G)= 0$ is stochastically greater 
than the conditional distribution given $e\_{\rm pos}(G) = 1$.
%\begin{figure}[h]
%    \centering
%       \includegraphics[width=0.6\textwidth]{degstd}
%    \caption{ }
%    \label{fig:}
%\end{figure}
\end{conj}

\begin{conj}[Stochastic Ordering of Number of Claws by $e\_{\rm pos}$]
For graphs in $\mathcal{G}_n$, the conditional distribution of number of claws given $e\_{\rm pos}(G) = 0$ is stochastically greater 
than the conditional distribution given $e\_{\rm pos}(G)= 1$.
\end{conj}

\firstlevel{Discussion and Evaluation of the Framework}

In this section, we address key points of discussion regarding the novelty, applicability, and findings of our proposed framework.

\secondlevel{On methodological novelty and significance}

Our framework builds upon the foundational work of Davies et al. \cite{Nature21}. A central question is whether our contribution constitutes a significant methodological advancement or merely an incremental adjustment. While the core architecture shares similarities, our contribution lies in a crucial conceptual and practical shift: moving from a framework designed to identify "if and only if" statements to one optimized for discovering sufficient conditions, framed as ''if ..., then ..." conjectures. This is not a trivial change, but a fundamental re-engineering of the learning objective to prioritize predictive precision. This prioritization steers the model away from finding necessary conditions and toward identifying high-precision, sufficient features.

The value of our framework is demonstrated by the new kinds of mathematical insights it facilitates. Applying the original Davies et al. framework to the problem of e-positivity did not yield results that suggested earlier known theorems. In contrast, our precision-focused framework reproduced evidence for a known theorem by Stanley and Stembridge—thereby validating its ability to align with mathematicians’ intuition.

\secondlevel{On the choice of problem and reproduction of results}

Our selection of e-positivity in algebraic combinatorics as a case study was intentional. This problem is well-studied, providing a rich landscape of established results against which to benchmark our novel framework. Reproducing known results, such as aspects of the Stanley-Stembridge theorem, is a critical validation step. It shows that the framework can independently and reliably recover mathematical truths from data, establishing its credibility before applying it to find conjectures.

\secondlevel{On the role of continuous vs. discrete invariants}

A key empirical finding of our work is the apparent importance of continuous graph invariants in predicting e-positivity. This claim is supported by our Saliency Map analysis, which consistently ranked continuous features among the most influential predictors. Specifically, the six highest-ranked features by this measure are continuous. In the context of our model, this indicates that small changes in these continuous invariants lead to larger changes in the predicted probability of $e$-positivity than changes in discrete features. We interpret this as evidence that the property of $e$-positivity, at least for the classes of graphs studied, is more sensitive to, and thus potentially more correlated with, continuous aspects of graph structure. This provides a novel, data-driven perspective on a property traditionally studied through discrete combinatorial lenses. Then, we find a way to formulate genuinely new conjectures and to relate $e$-positivity to continuous invariants in a more mathematical manner, rather than relying solely on data-driven perspectives.

\firstlevel{Acknowledgements} We thank Stephanie van Willigenburg for many informative and thought-provoking conversations. Both authors were supported in part by the Provincial Nature Science Foundation of Shandong, Project No. ZR2024QA026 and the Fundamental Research Funds for the Central Universities.

\firstlevel{Data availability statement}\label{code}
The notebooks and the data are available at\\  \href{https://github.com/AIMath-Lab/SufficientConditions}{https://github.com/AIMath-Lab/SufficientConditions}.

%
%
%  \noindent   Farid Aliniaeifard ({\href{mailto:farid@sdu.edu.cn}{farid@sdu.edu.cn})  \\
%{\textit{Research Center for Mathematics and Interdisciplinary Sciences, Shandong University}}\\
%{\textit{Frontiers Science Center for Nonlinear Expectations, Ministry of Education, Qingdao, Shandong,
%266237, P. R. China}}\\
%\\
%Shu Xiao Li ({\href{mailto:lishuxiao@sdu.edu.cn}{lishuxiao@sdu.edu.cn})  \\
%{\textit{Research Center for Mathematics and Interdisciplinary Sciences, Shandong University}}\\
%{\textit{Frontiers Science Center for Nonlinear Expectations, Ministry of Education, Qingdao, Shandong,
%266237, P. R. China}}
%

\end{document}